%% file: main.tex
\documentclass[runningheads]{llncs}

\usepackage{amsmath}
\usepackage{graphicx}
\DeclareMathOperator*{\argmin}{arg\,min}

\begin{document}

\title{Minimizing Entropy to Discover Good Solutions to Recurrent Mixed Integer Programs}

\titlerunning{Minimizing Entropy to Discover Good Solutions to Recurrent MIPs}

%
%
\author{Charly Robinson La Rocca\inst{1}
\and
Emma Frejinger\inst{1} 
\and
Jean-François Cordeau \inst{2}
}
\authorrunning{C. Robinson et al.}
%
\institute{Université de Montréal, Montréal, Canada 
\email{\{charly.robinson.la.rocca,emma.frejinger\}@umontreal.ca}
\and
HEC Montréal, Montréal, Canada \\
\email{jean-francois.cordeau@hec.ca}}
\maketitle              
\begin{abstract}
Current state-of-the-art solvers for mixed-integer programming (MIP) problems are designed to perform well on a wide range of problems. However, for many real-world use cases, problem instances come from a narrow distribution.  This has motivated the development of specialized methods that can exploit the information in historical datasets to guide the design of heuristics. Recent works have shown that machine learning (ML) can be integrated with an MIP solver to inject domain knowledge and efficiently close the optimality gap.  This hybridization is usually done with deep learning (DL), which requires a large dataset and extensive hyperparameter tuning to perform well. This paper proposes an online heuristic that uses the notion of entropy to efficiently build a model with minimal training data and tuning.  We test our method on the locomotive assignment problem (LAP), a recurring real-world problem that is challenging to solve at scale. Experimental results show a speed up of an order of magnitude compared to a general purpose solver (CPLEX) with a relative gap of less than 2\%. We also observe that for some instances our method can discover better solutions than CPLEX within the time limit.

\end{abstract}

\input{introduction}

\input{related_works}
\input{methodology}
\input{experiments}

\input{conclusion}

\newpage 
\bibliographystyle{splncs04}
\bibliography{mybib}

\end{document}

%% file: introduction.tex
\section{Introduction}


Machine learning (ML) has made major leaps forward with the maturation of specialized hardware \cite{wang2019benchmarking} and software \cite{paszke2019pytorch} stacks.  ML enables practitioners to automatically find useful patterns from a dataset without the extensive domain knowledge often required to make progress on a given problem. In the context of our research, we are interested in the exploitation of learning algorithms to efficiently solve large combinatorial optimization (CO) problems. The idea to leverage the knowledge from historical data is not new. Many strategies have emerged in recent years to take advantage of learning methods and improve upon standard solvers.  \cite{Bengio_2020} provide a tour d'horizon on ML for combinatorial optimization and give an overview on the different ways in which they can be integrated.  The first idea uses ML to directly predict the solution of the problem instance from the input features. One example of this is the application of the pointer network developed by \cite{vinyals2015pointer} to learn to predict the optimal permutation for the travelling salesman problem (TSP). The second approach is to learn to optimize the hyperparameters of the CO algorithm. Third, ML can also be used alongside the optimization algorithm. In practice, this is often done via the integration of learning in the branch-and-bound (B\&B) framework \cite{lodi2017learning}. Our approach belongs to the third class of hybrid methods.\\

This paper proposes an approach that focuses on simplicity and efficiency. Unlike most recent papers on the subject, we avoid using any form of deep learning (DL). Instead, we exploit the well-known notion of entropy to make decisions that \emph{minimize risk} when heuristically adding cuts to the MIP. We favour classic ML models to fit our data because of their robustness; they require significantly fewer data points and tuning than a deep neural network (DNN).  The model is updated \emph{online} during the B\&B and its predictions are used to guide the decision of adding cuts to the MIP. The goal is to remove from the search space the solutions that are the least likely to be optimal. We apply this strategy to the locomotive assignment problem (LAP), a problem that is not only challenging to solve but also \emph{recurring} i.e., similar problem instances are solved repeatedly over time. Hence, there is a substantial economic incentive to find a heuristic that can solve the LAP quickly.  We use the Canadian National Railway Company (CN) as the subject of our case study for the LAP. \\

This document first presents some related work on the integration of ML for CO problems. This is followed by an overview on how we use the notion of entropy to reduce the computational load necessary to solve the LAP. Finally,  we include some preliminary results for the speed up and the gap relative to our baseline CPLEX.

%% file: related_works.tex
\section{Related works}
One defining property of NP-hard CO problems is the fact that there are no algorithms with guarantees that can solve them in polynomial time unless $\text{P}= \text{NP}$. In practice, this means that they may be intractable when the size of the instance is large. This is why efficient heuristics are sought after. In this case, we pay the penalty of no optimality guarantees to find solutions in a reasonable amount of time. The design of heuristics requires substantial problem-specific research and trial-and-error on the part of practitioners. This constitutes the main motivation to use ML, since it can extract valuable patterns from data without domain knowledge about the given problem. \\

The first question that arises when integrating ML into CO is how to represent the input instance. CO problems may comprise both continuous and integer variables where the relationship is often modelled with graph data structures.  Unlike images, audio, and text, which have a clear grid structure, graphs have irregular structures,  making it challenging  to  generalize  some  of  the  basic mathematical operations to them. The dominant and fast-growing paradigm for DL with graph data is the GNN formalism \cite{hamilton2020graph}. There have been many successful attempts at using GNNs to solve CO problems \cite{peng2020graph,vessel}. We have also identified recent papers on GNNs that focus on how to formally improve their expressiveness, since there is a limited understanding of their representational properties and limitations \cite{xu2018powerful,alon2021on}. There is a bottleneck issue where all the information from the related nodes is compressed into a fixed-length vector and therefore GNNs fail to propagate messages flowing from distant nodes.  By way of contrast, our method is able to learn from a few features that are designed to be tightly correlated with the variables that we want to predict. This allows us to use a non-parametric model that can be trained efficiently compared to DNNs.  \\

During the B\&B, there are many different types of decisions that have to be made at each node in the tree. The two most common ones are associated with the branching scheme: variable selection (which fractional variable to branch on) and node selection (which of the currently open nodes to process next). The popular framework used to integrate ML within the B\&B is RL \cite{gasse2019exact,cappart2020combining,dai2017learning,barrett2020exploratory}. This is expected because there is a natural fit between RL and the iterative algorithm used to solve CO problems. We can model the branching decisions as actions made by the agent.  The state of the environment can be modelled as the input instance with the current position in the B\&B tree.  Nevertheless, RL has limitations that make it challenging to implement. The main concern comes from the design of the reward function. Learning in the RL framework is mainly driven by the reward signal which associates a reward value to each transition.  Some use a reward of -1 at each step to motivate exploration but that reward function contains no guiding signal \cite{SeaPearl}. Moreover, training an agent using RL may be unstable and tricks such as replay memory and gradient clipping are often necessary to achieve good performance. For these reasons, instead of using RL to train our policy, we design it around the idea of risk minimization. We use the entropy of the variable to measure how risky it is to add a lazy constraint to it. The simplicity of our method makes the training process trivial, significantly improving its robustness compared to RL.  \\

In summary, we noticed that GNNs and RL are popular tools used to learn to solve CO problems. However, they require extensive training and tuning to yield acceptable performance.  In contrast, our method targets the idea of entropy minimization to directly capture the relevant information from historical data. As a consequence, we can take advantage of simpler models that are less prone to over-fitting. Furthermore, the online nature of the model means that most of the learning happens during the B\&B of the test instance. Those design choices creates the opportunity for an accurate model that requires minimal training data and manual tuning.


%% file: methodology.tex
\section{Methodology}\label{sec:methodology}
\newcommand{\V}{\mathcal{V}}

We start with a general MIP formulation which is described using the notation $ \argmin_{\mathbf{x}} \{ c^T \mathbf{x}\,  |\,  A\mathbf{x} \leq b, \, \mathbf{x} \in \mathcal{B} \cup \mathcal{Q} \cup \mathcal{P}  \} $,  where $\mathbf{x} $ is the vector of decision variables that is partitioned into  $\mathcal{B}$, $\mathcal{Q}$ and $\mathcal{P}$, the index sets of binary, integer and continuous variables, respectively. To be compatible with our methodology, the MIP needs to have at minimum the following constraints 
\begin{align}
    & \sum_{k \in K} x_v^k = 1  &\forall v \in \V,  \label{c_class}
\end{align}
where the binary variable $x_v^k \in \mathcal{B}$ is equal to 1 when we assign the class $k \in K$ to the object $v \in \V$. The constraints (\ref{c_class}) ensure that only one class is possible per object. We use the notation $k_{vt}$ to refer to the class of $v$ at iteration $t$ during the B\&B.  We characterize our methodology as general because it is often trivial to transform a MIP to fit this formulation. In theory, one can choose to convert an integer variable into a set of binary variables without changing the problem definition. The idea is to use a ML model to predict the class of some objects and let the MIP solver optimize all the other variables. If we fix the binary variables to be equal to their optimal value, we expect that the time required to solve the MIP will be greatly reduced. \\


\textbf{Representation embedding.} The current method uses an \emph{online} learning algorithm that updates its belief dynamically during the B\&B as new data becomes available. This is a significant change compared to the standard supervised learning approach where the dataset is generated offline \cite{khalil,Gupta20hybrid}. In the online paradigm, the feature vector needs to be updated at every iteration. Let  $\phi_t(v)$  be the feature vector associated with variable $v$ in the MIP at iteration $t$ in the B\&B.  To avoid performance concerns, the complexity of the function that computes  $\phi_t(v)$  given $\phi_{t-1}(v)$  should be $O(1)$. For this reason, $\phi_t(v)$ is calculated using statistics that can be easily updated such as the mean, the variance and the extrema (min and max). These statistics are evaluated on the solutions collected during B\&B up until the current iteration $t$: $K_{vt} =\{ k_{v1},k_{v2},...,k_{vt} \}$. Formally, the feature vector $\phi_t(v)$ is calculated using the expression 
\begin{align}
    \phi_t(v) = \left[\text{mean}(K_{vt}), \text{var}(K_{vt}), \text{max}(K_{vt}), \text{min}(K_{vt})\right].
\end{align}
The next step is to discuss the corresponding label $s_v$ associated with $v$. The label is binary and informs us about the \emph{stability} of the variable:  $s_v=1$ when the variable $v$ is stable and $s_v=0$ otherwise.  The idea of using the concept of stability as a label was inspired by \cite{ding2020accelerating}. We build on this idea by using the notion of entropy $H$ to infer the stability \cite{shannon1948mathematical}. It is defined as $H(Z) = - \sum_i P(z_i) \log P(z_i)$ where $P(z_i)$ is the probability associated with the outcome $z_i$ for the the random variable $Z$.  We convert the entropy into a binary variable using the following function 


\begin{equation}
    s_v=
    \begin{cases}
      0, & \text{if}\ H(K_v) > m \\
      1, & \text{otherwise,}
    \end{cases}
\end{equation}
where $H(K_v)$ is the entropy of the visited solutions and $m$ is the median of the entropy for each variable in the instance. To generate a balanced training set, we use the median because it reliably creates a 50/50 split between stable and unstable variables.\\

\textbf{Machine learning model.} The notable attribute of our model is its online nature. Also,  the model is potentially called during each callback of the solver, which incentivizes a lean design.  We decided to use the Fast Random Forest (FRF) model from the OnlineStats.jl library written in \emph{Julia} \cite{Day2020}. This ensemble model has fast processing speed that can be updated in constant time. Furthermore, the advantage of having the whole software stack written in \emph{Julia} significantly reduces the overhead that we usually have with libraries in Python and C++.    \\

\textbf{Policies.} A policy decides which action to take given the predictions generated with the ML model. In our application, the action space is defined by the interface of the general purpose solver CPLEX. We can name a few different types of actions such as user cuts, lazy constraints and branching rules. We selected lazy constraints as our main tool to interact with the solver because it limits the solution directly, which greatly reduces the computational time when it is done appropriately. \\

We propose some basic rules and assumptions to significantly simplify the design space of policies. First, we assume that lazy constraints can only be added after an incumbent solution is found. Second, the policy should be risk averse to avoid removing the optimal solution from the search space.  Third, we pose that each lazy constraint should have the following form: $x_v^{k} = 1$ where $k$ is the predicted class i.e., $k = \tilde{k}_v$. We use these assumptions to design two different types of policies: a scoring policy that uses the entropy directly and a threshold policy that uses the predictions of the FRF model to select the variables to fix. \\


The \emph{scoring policy} (SP) uses a scoring system to measure the risk associated with fixing each variable in the MIP. At each incumbent, SP computes the score $\sigma_v$ for each variable $v$ as follows: $ \sigma_v = - H(K_{vt})$ where $H(K_{vt})$ is the entropy of the solutions visited up until the current iteration $t$. We use the negative because we want to maximize the score. In other words, we are interested in the variables that are the least likely to change. Once the scores are computed, we select the top $n$ from the list and add a lazy constraint to each one of them, where $\tilde{k}_v$ is the most common class in $K_{vt}$. For example, SP($n=5$) will select five variables to fix at each incumbent using the entropy based scoring system. \\

The \emph{threshold policy} (TP) is a more aggressive version of SP because every constraint is added after the first incumbent. TP uses a threshold $\tau$ to filter which variable should be fixed based on the prediction history of the ML model. A prediction $\tilde{s}_{vt}$ is made at every iteration $t$ in the B\&B and it is stored in what we call the prediction history $\tilde{S}_v = \{\tilde{s}_{v1}, \tilde{s}_{v2}, ..., \tilde{s}_{vt} \}$. If the model systematically predicts that the variable is stable, the policy will add a constraint to it. More formally, the set of selected variables is $\{ v | \text{mean}(\tilde{S}_v) > \tau \}$. For example, TP($\tau = 0.5$) will select all the variables for which the ML model predicts $\tilde{s}_v=1$ more than 50\% of the time. 

%% file: experiments.tex
\section{Experiments}

We apply our methodology to the LAP using the historical schedules of CN. This problem is difficult to solve without partitioning it into daily sub-problems, because it contains a large number of binary variables. It follows daily and weekly cyclical patterns which are useful for learning. The goal of the LAP is to find the cheapest way to send a certain number of locomotives through a network such that all capacity and demand constraints are satisfied. The decisions $ x_v^k$ represent the configuration of locomotives $k$ assigned to each train arc $v$ in the instance. Other variables are required to model flow conservation constraints for each node in the network. We refer to \cite{camilo} for the details on the MIP formulation of the LAP. \\

\textbf{Setup. } 
We define three instance sizes: small (100 to 250 trains), medium (250 to 450 trains) and large (450 to 800 trains). For each size, we generate 10 consecutive weekly instances that are solved with a time limit of 30 minutes. We compare the baseline CPLEX with three policies: $SP(n=1)$, $SP(n=5)$ and $TP(\tau =0.5)$. Our software stack is fully written in \emph{Julia} and we use the \emph{JuMP} library to interact with the solver \cite{jump}.  We run our experiments on a t2.xlarge instance on Amazon Web Services (AWS) with 4 vCPU and 16 GB of memory. The FRF model is trained on a single small instance to learn to predict the stability of each decision variable. The training procedure takes 4.54 seconds to complete.  \\

\textbf{Metrics. } We consider the work of \cite{berthold2013measuring}, which is focused on this specific issue of measuring the impact of primal heuristics. They introduce a performance measure that takes into account the whole solution process by computing the integral of the primal gap over time. They define the primal integral $PI$ of a run as $ PI := \sum_i p(t_{i-1})\cdot (t_i - t_{i-1})$, where $t_i $ is the point in time when the $i$th new incumbent solution is found and $p(t_i)$ is the primal gap associated with it. The smaller $PI$ is, the better the expected quality of the solution if we stop the solver at an arbitrary point in time. To get a relative measure of performance, we propose the \emph{primal integral ratio} (PIR) which is the ratio between the primal integral of the heuristic and the baseline. To provide a good intuitive understanding for how much faster the heuristic is, we also compute a custom version of the speed up. Most instances reach the time limit of 30 minutes, therefore we compute its value at every node in the B\&B. By leveraging a linear interpolation, we can create a function that returns the speed up for a given gap. From this, we can evaluate an upper bound on the speed up and its associated gap; we call them \emph{best speed up} and \emph{gap at best speed up} respectively. We should mention that the gap is measured relative to the best feasible solution found by CPLEX within the time limit.  \\

\begin{figure}[!ht]
    \begin{center}
        \begin{minipage}{0.45\textwidth}
        \centering
            \scalebox{1.1}{\includegraphics[width=\textwidth]{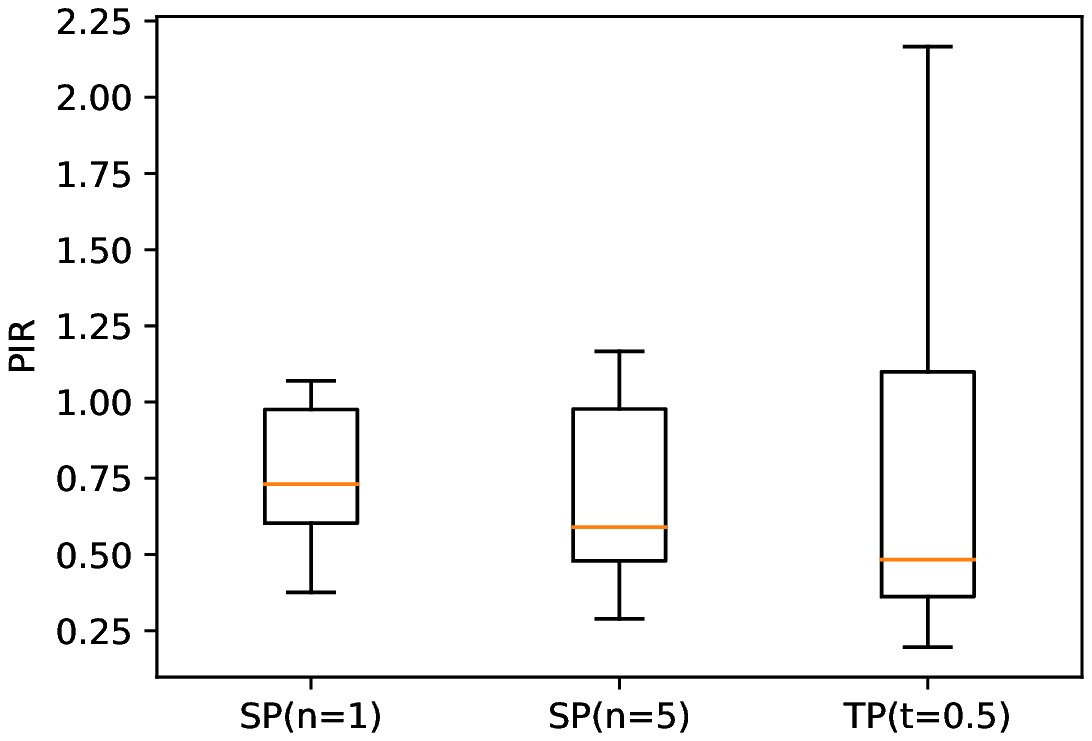}}
        \caption{PIR for different policies} \label{fig:pir_policy}    
        \end{minipage}
        \hfill
        \begin{minipage}{0.45\textwidth}
        
        \centering
            \scalebox{1.1}{\includegraphics[width=\textwidth]{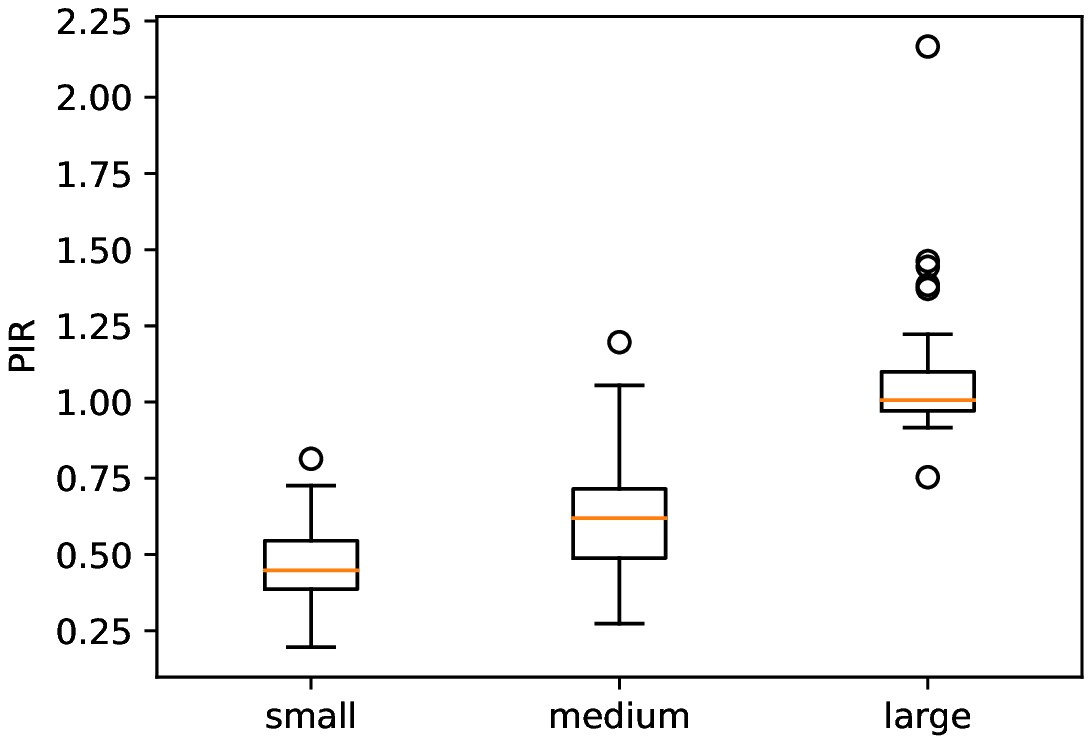}}        
            \caption{PIR for different instance sizes} \label{fig:pir_instance_size}
    \end{minipage}
    \end{center}
\end{figure}

\begin{figure}[!ht]
    \begin{center}
        \begin{minipage}{0.45\textwidth}
        \centering
            \scalebox{1.1}{\includegraphics[width=\textwidth]{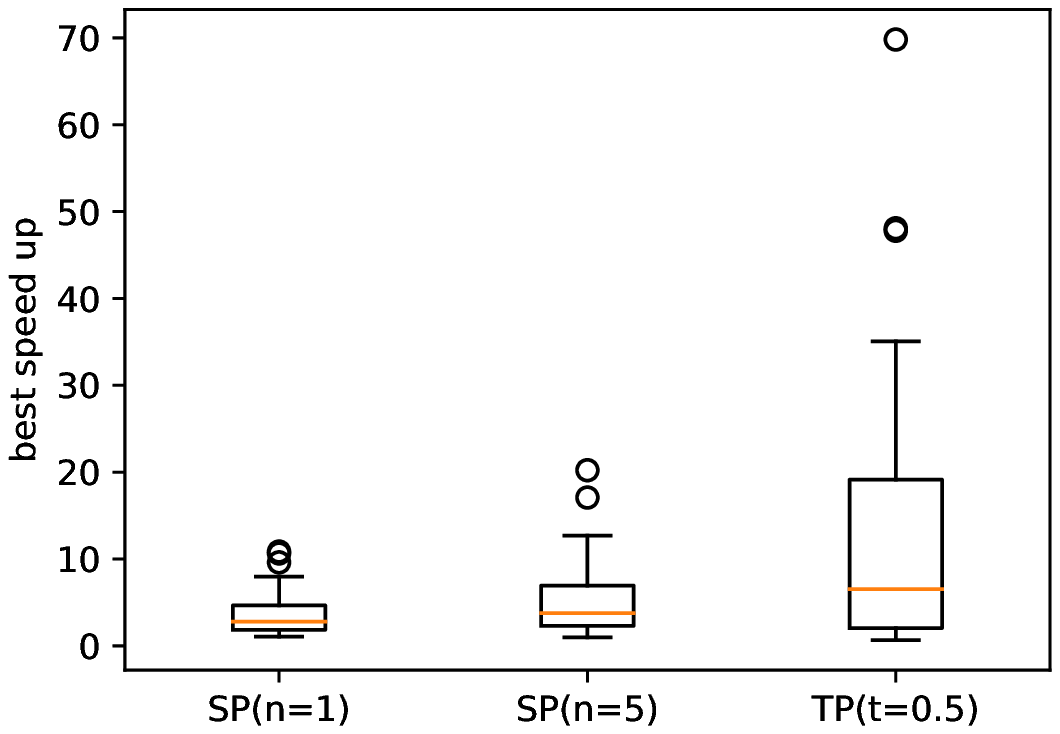}}
        \caption{Best speed up for different policies} \label{fig:speed_up}    
        \end{minipage}
        \hfill
        \begin{minipage}{0.45\textwidth}
        
        \centering
            \scalebox{1.1}{\includegraphics[width=\textwidth]{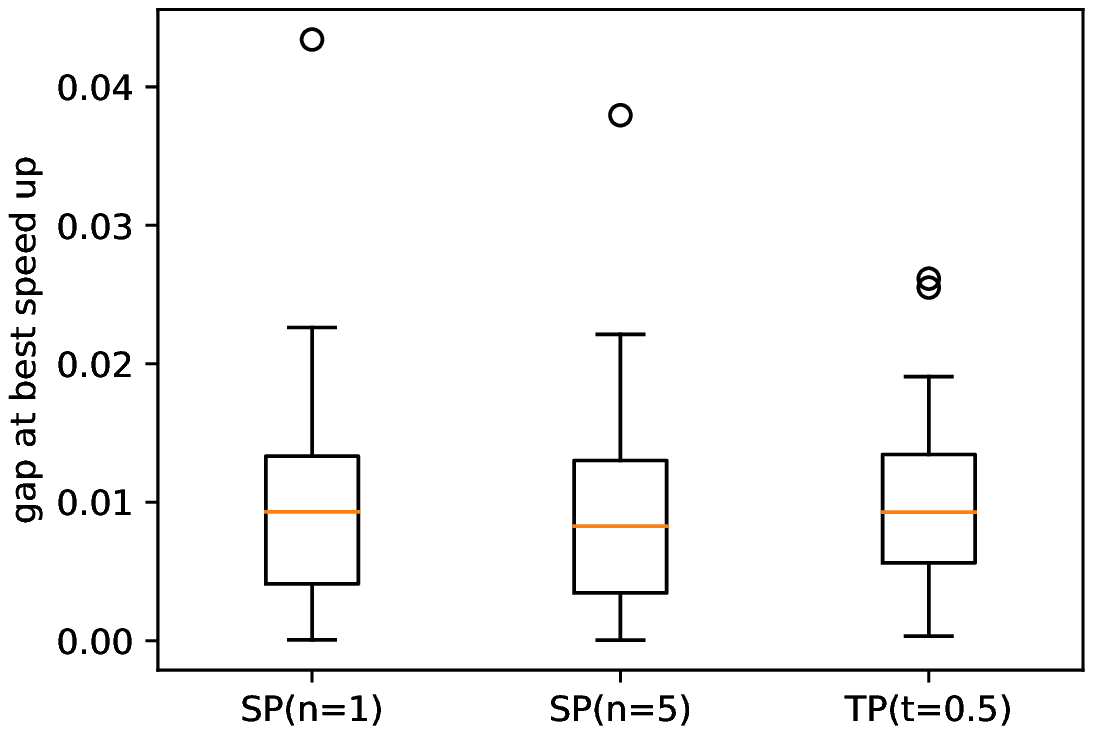}}        
            \caption{Gap at best speed up for different policies } \label{fig:gap}
    \end{minipage}
    \end{center}
\end{figure}

\textbf{Results. } For the metric PIR we group the results by policy and instance size. The body on the box plot illustrate the 25th and 75th percentiles of the distribution. Fig.  \ref{fig:pir_policy} shows that the 75th pencentile of PIR is 0.969, 0.977 and 1.17 for the policies $SP(n=1)$, $SP(n=5)$ and $TP(\tau =0.5)$ respectively.  Since for most instances the PIR is below 1, it means that we converge faster than CPLEX in general. In Fig.  \ref{fig:pir_instance_size}, we have a similar pattern for small and medium instances as their box plot body is well below 1. There is a noticeable performance drop for large instances (25th and 75th percentiles are 0.972 and 1.15 for the PIR) which is expected given that we train our model on a small instance. Surprisingly, we measure a negative relative gap for 66.7\% (20/30) of large instances, which suggests that our method can also be useful to find a better feasible solution within a limited time frame. Figs. \ref{fig:speed_up}-\ref{fig:gap} illustrate the speed up and gap respectively. TP is the most aggressive one and this is noticeable here, as we observe that the limit for outliers is 38x.  The 75th percentile of the gap is below 2\% for all policies, which is reasonable for users that value speed over quality. The good performance of our method can be explained by the high \emph{action accuracy} i.e., the ratio of lazy cuts that preserve the optimal solution. We measure an action accuracy per instance of 94.4\% on average and 80.4\% in the worst case.  We ought to mention that we experimented with supervised learning methods that directly predict the solution of instances using the problem graph as input. This approach had a substantially worst accuracy (between 50\% and 75\%) and, as a result, the solution quality was unsatisfactory. 

%% file: conclusion.tex
\section{Perspectives and future work}
There are many open questions with respect to the integration of ML into CO solvers. There is the important question of how to represent the input data of the problem instance. According to our reading of the literature, GNNs or more generally DNNs are the most popular tools used for representation learning.   Given enough data, they can learn a high quality embedding of the input space without requirements for prior domain knowledge from practitioners.  However, they come with practical challenges that impinge upon the research process. DNNs are more costly and data intensive to train than the model considered in this work.\\ 

The complexity overhead of DL motivated our design of a new method that decides how to add cuts to a MIP via entropy minimization. Given the online aspect of the model, most of the learning happens at test time and, therefore, it requires a negligible amount of resources to train in terms of both data and time. A method like ours, that takes less than 5 seconds to train, allows us to iterate faster on the core mechanic of the heuristic. This strategy ended up being successful, as we were able to find feasible solutions for the LAP with a relative gap below 2\% up to 10 times faster than CPLEX. There are still many paths to explore to validate the methodology. We plan to apply it to other similar CO problems such as the \emph{load planning problem} \cite{lpp}. We identified the containers-to-car matching problem as a potential entry point for our learned heuristic.   Once it is validated on multiple problems, the next step is to tackle the scaling challenges. Future work will be dedicated to issues occurring when solving even larger instances.  \\